\documentclass[11pt]{article}
\input epsf
\usepackage{amssymb}
\usepackage{amsmath}
\usepackage{amsfonts}
\input amssym.def
\usepackage{amssymb}
\begin{document}
\title{\bf The well-posedness of Cauchy problem for
dissipative modified Korteweg de Vries equations
\thanks{The work was partly supported by NNSF of China (No.10371080)
and China Postdoctoral Science Foundation (No.\,2005038327).}}
\author{ Wengu Chen$^1$, Junfeng Li$^{1,\,2}$ and Changxing Miao$^1$ \\$^1$Institute of Applied Physics
and Computational Mathematics\\P.O.Box 8009, Beijing 100088,
China\\$^2$College of Mathematics, Beijing Normal University,
\\Beijing 100875, China}

\date{E-mail:\,chenwg@iapcm.ac.cn,\, junfli@yahoo.com.cn,\, miao$_-$changxing@iapcm.ac.cn}

\maketitle

{\bf Abstract.}\, In this paper we consider some dissipative
versions of the modified Korteweg de Vries equation
$u_t+u_{xxx}+|D_x|^{\alpha}u+u^2u_x=0$ with $0<\alpha\leq 3$. We
prove some well-posedness results on the associated Cauchy problem
in the Sobolev spaces $H^s({\Bbb R})$ for $s>1/4-\alpha/4$ on the
basis of the $[k;\,Z]-$multiplier norm estimate obtained by Tao in
\cite{Tao} for KdV
equation. \\
{\bf 2000 Mathematics Subject Classification:}\, 35Q53, 35A07.

\section{Introduction} The $X^{s,\,b}$ spaces, as used by
Beals, Bourgain, Kenig-Ponce-Vega, Klainerman-Machedon and others,
are fundamental tools to study the low-regularity behavior of
nonlinear dispersive equations. It is of particular interest to
obtain bilinear or multilinear estimates involving these spaces.
By Plancherel's theorem and duality, these estimates reduce to
estimating a weighted convolution integral in terms of the $L^2$
norms of the component functions. In \cite{Tao}, Tao
systematically studied the weighted convolution estimates on
$L^2$. As a consequence, Tao obtained sharp bilinear estimates for
the KdV, wave, and Schr\"odinger $X^{s,\,b}$ spaces. These
estimates are usable for many applications, and Tao in his paper
presented some selected applications of these estimates to prove
multilinear estimates and local well-posedness results. In this
paper, we will get some bilinear estimate on some new type of
$X^{s,\,b}$ spaces adapted to the dissipative version of the
modified Korteweg de Vries equation and obtain new local
well-posedness of associated Cauchy problem. This new type of
Bourgain spaces was first introduced by Molinet and Ribaud in
\cite{MR1}. We defer to introduce Tao's estimates here.

Recently, Molinet and Ribaud consider the Cauchy problem
associated with dissipative Korteweg de Vries equations
\begin{eqnarray}
\left\{\begin{array}{ll} u_t+u_{xxx}+|D_x|^\alpha u+uu_x=0,\,t\in
{\Bbb
R}_+,\,x\in{\Bbb R},\\
u(0)=\varphi,\end{array}\right.\end{eqnarray} where $|D_x|^\alpha$
denotes the Fourier multiplier operator with symbol $|\xi|^\alpha$.
When $\alpha=1/2$, equation (1) models the evolution of the free
surface for shallow water waves damped by viscosity. When
$\alpha=2$, equation (1) is the so-called Korteweg de Vries Burgers
equation which described the propagation of small amplitude long
waves in some nonlinear dispersive media when dissipative effects
occur (see\cite {OS}).

In \cite{MR1}, Molinet and Ribaud proved the global well-posedness
of (1) for data in $H^s({\Bbb R}),\,s>-\frac 34$ for all $\alpha>0$.
Especially when $\alpha=2$, they proved the global well-posedness of
Korteweg de Vries Burgers equation for data in $H^s({\Bbb R}),\,s>
-\frac 34-\frac 1{24}$. The surprising part of this result is that
the indice $s=-\frac 34-\frac 1{24}$ is lower than the best known
indice $s=- 3/4$ obtained by Kenig, Ponce and Vega in \cite{KPV1}
for the KdV equation and lower than the indice $s=- 1/2$ of the
critical Sobolev space for the dissipative Burgers equation
$u_t-u_{xx}+uu_x=0$ (see \cite{B}, \cite{Dix2}). In \cite{MR2},
Molinet and Ribaud improved this result by introducing a new
Bourgain type space and working in the space. They showed that
Korteweg de Vries Burgers equation is globally well-posed in
$H^s({\Bbb R})$ for $s>-1$ and in some sense ill-posed in $H^s({\Bbb
R})$ for $s<-1$.

The main purpose of this paper is to consider the Cauchy problem
for the following dissipative versions of the mKdV equation on the
real line
\begin{eqnarray}\label{2}
\left\{\begin{array}{ll} u_t+u_{xxx}+|D_x|^{\alpha}u+u^2u_x=0,\,t\in{\Bbb R}_+\\
 u(0)=\varphi \end{array}\right.\end{eqnarray}
 with $\alpha\in (0,\,3]$. And we prove local well-posedness results
 on the associated
Cauchy problem in the Sobolev spaces $H^s({\Bbb R})$ for
$s>1/4-\alpha/4$ on the basis of the $[k;\,Z]-$multiplier norm
estimate obtained by Tao in \cite{Tao} for KdV equation. We also get
the global well-posedness for $1<\alpha\leq 3$ when
$s>\frac{1}{4}-\frac{\alpha}{4}$. It is worth pointing out that the
method used to obtain the multilinear estimates by Molinet and
Ribaud in \cite{MR2} does not work when $\alpha$ is small enough. If
we consider the case $0<\alpha\leq 1$, we can only get that problem
(2) is local well-posed for $s>\frac{1}{2}-\frac{\alpha}{2}$ by
running the approach of \cite{MR2}. The basic ideas of processing
the multilinear estimates on $X^{s,\,b}$ spaces are using
Cauchy-Schwarz inequality and the algebraic smoothing relation in
\cite{MR2}( resonance identity in \cite{Tao}). In \cite{Tao}, Tao
utilizes dyadic decomposition and orthogonality before resorting to
Cauchy-Schwarz. The advantages of dyadic decomposition and
orthogonality lead to a better estimate when the algebraic smooth
relation brings little benefit.

It is clear we get a better local well-posedness for the
dissipative version of the mKdV equation than the local
well-posedness of mKdV. This explicates that the dissipative term
in (\ref{2}) plays a key role for the low regularity of the
equation. We are interested in finding out the smooth effect of
the dissipative term. From the scaling analysis, the critical
index of the mKdV equation is $-\frac{1}{2}$, while the critical
index of the equation with only dissipative term is
$1-\frac{\alpha}{2}$. Thus we conjecture that the mean value
$\frac{1}{4}-\frac{\alpha}{4}$ of the two critical indices is the
critical index for our dissipative version of the mKdV. For
instance, we consider the Korteweg de Vries Burgers equation, $-1$
is exactly the mean value of $-\frac{3}{2}$ the critical scaling
index of KdV and $-\frac{1}{2}$ the critical scaling index of
Burgers equation. Unfortunately, we can not show the ill-posedness
below the index $\frac{1}{4}-\frac{\alpha}{4}$. Recently, the
first two authors of this paper \cite{CL} showed that when
$\alpha=2$, (\ref{2}) is ill-posedness for $s<-\frac{1}{2}$ in
some sense. A similar argument can also get the same ill-posedenss
for our problem.

\subsection{Notations}
For a Banach space $X$, we denote by $\|\cdot\|_X$ the norm in
$X$. We will use the Sobolev spaces $H^s({\Bbb R})$ and their
homogeneous versions $\dot{H}^s({\Bbb R})$ equipped with the norms
$\|u\|_{H^s}=\|(1-\Delta)^{s/2}u\|_{L^2}$ and
$\|u\|_{\dot{H}^s}=\||D|^{s/2}u\|_{L^2}$. Recall that for
$\lambda>0$,
\begin{eqnarray}
\|f(\lambda t)\|_{H^s}\leq
(\lambda^{-1/2}+\lambda^{s-1/2})\|f(t)\|_{H^s},\,\|f(\lambda
t)\|_{\dot{H}^s}\sim \lambda^{s-1/2}\|f(t)\|_{\dot{H}^s}.
\end{eqnarray}
We also consider the corresponding space-time Sobolev spaces
$H^{s,\,b}_{x,\,t}$ endowed with the norm
\begin{eqnarray}
\|u\|^2_{H^{s,\,b}}=\int_{{\Bbb
R}^2}<\xi>^{2s}<\tau>^{2b}|\hat{u}(\xi,\,\tau)|^2d\xi d\tau,
\end{eqnarray}
where $<\cdot>=(1+|\cdot|^2)^{1/2}$.

Let $U(\cdot)$ be the unitary group in $H^s({\Bbb R}),\,s\in{\Bbb
R}$, which defined the free evolution of the KdV equation, i.\,e.
\begin{eqnarray}
U(t)=\exp\Big(itP(D_x)\Big),
\end{eqnarray}
where $P(D_x)$ is the Fourier multiplier with symbol
$P(\xi)=\xi^3$. Since the linear symbol of equation (2) is
$i(\tau-\xi^3)+|\xi|^\alpha$, by analogy with the spaces
introduced by Bourgain in \cite{Bourgain} for purely dispersive
equations
 and Molinet and Ribaud for KdV-Burgers equation,
we define the function space $X^{s,\,b}_\alpha$ endowed with the
norm
\begin{eqnarray}
\|u\|_{X^{s,\,b}_\alpha}=\|<i(\tau-\xi^3)+|\xi|^\alpha>^b<\xi>^s\hat{u}(\xi,\,\tau)\|_{L^2({\Bbb
R}^2)}.
\end{eqnarray}
So that
\begin{eqnarray}
\|u\|_{X^{s,\,b}_\alpha}\sim\|<|\tau-\xi^3|+|\xi|^\alpha>^b<\xi>^s\hat{u}(\xi,\,\tau)\|_{L^2({\Bbb
R}^2)}.
\end{eqnarray}

Note that since ${\cal F}(U(-t)u)(\xi,\,\tau)={\cal
F}(u)(\xi,\,\tau+\xi^3)$, one can re-express the norm of
$X^{s,\,b}_\alpha$ as
\begin{eqnarray}
\|u\|_{X^{s,\,b}_\alpha}&=&\|<i\tau+|\xi|^\alpha>^b<\xi>^s\hat{u}(\xi,\,\tau+\xi^3)\|_{L^2({\Bbb
R}^2)}\nonumber\\
&=&\|<i\tau+|\xi|^\alpha>^b<\xi>^s{\cal
F}(U(-t)u)(\xi,\,\tau)\|_{L^2({\Bbb R}^2)}.
\end{eqnarray}
For $T\geq 0$, we consider the localized spaces
$X^{s,\,b}_{\alpha,\,T}$ endowed with the norm
$$\|u\|_{X^{s,\,b}_{\alpha,\,T}}=\inf_{w\in
X^{s,\,b}_\alpha}\{\|w\|_{X^{s,\,b}_\alpha},\,w(t)=u(t)\,{\rm
on}\,[0,\,T]\}.$$ Finally we denote by $W(\cdot)$ the semigroup
associated with the free evolution of the equation (2), i.\,e.
$$\forall t\geq 0,\,{\cal
F}_x(W(t)\varphi)(\xi)=\exp[-t|\xi|^\alpha+it\xi^3]\hat{\varphi}(\xi),$$
and we extend $W(\cdot)$ to a linear operator defined on the whole
real axis by setting
\begin{eqnarray}
\forall t\in{\Bbb R},\,{\cal
F}_x(W(t)\varphi)(\xi)=\exp[-|t||\xi|^\alpha+it\xi^3]\hat{\varphi}(\xi).
\end{eqnarray}

\subsection{Main results}
We will mainly work on the integral formulation of (2), i.\,e.
\begin{eqnarray}
u(t)=W(t)\varphi-\frac
13\int_0^tW(t-t')\partial_x(u^3(t'))dt',\,t\geq 0.
\end{eqnarray}
Actually, to prove the local existence result, we shall apply a
fixed point argument to the following truncated version of (10)
\begin{eqnarray}
u(t)=\psi(t)\Big[W(t)\varphi-\frac {\chi_{{\Bbb R}_+}(t)}
3\int_0^tW(t-t')\partial_x(\psi_T^3(t')u^3(t'))dt'\Big],
\end{eqnarray}
where $t\in {\Bbb R}$ and, in the sequel of this paper, $\psi$ is
a time cut-off function satisfying
$$\psi\in C_0^\infty({\Bbb R}),\,{\rm supp}\,\psi\subset [-2,\,2],\,
\psi\equiv 1\,{\rm on}\,[-1,\,1],$$ and
$\psi_T(\cdot)=\psi(\cdot/T)$. Indeed, if $u$ solves (11) then $u$
is a solution of (10) on $[0,\,T],\,T<1$.

Let us first state our global well-posedness result on the real
line. In this paper we restrict us with $0<\alpha\leq 3$.
\newtheorem{theorem}{Theorem}
\begin{theorem} Let $\varphi\in H^s({\Bbb R}),\,s>1/4-\alpha/4$. Then there exist
some $T>0$ and a unique solution $u$ of (10) in
\begin{eqnarray}
Z_T=C([0,\,T],\,H^s)\cap X_{\alpha,\,T}^{s,\,1/2}.
\end{eqnarray}
Moreover the map $\varphi\mapsto u$ is smooth from $H^s({\Bbb R})$
to $Z_T$.
\end{theorem}

If the dissipative term is strong enough, we can also get the global
well-posedness.
\begin{theorem}Let $\varphi\in H^s({\Bbb R}),\,s>1/4-\alpha/4$ and $1<\alpha\leq3$.
The existence time of the solution of (10) can be extended to
infinity. And $u$ belongs to $C((0,\,+\infty),\,H^\infty({\Bbb
R}))$.
\end{theorem}

{\bf Remark}: {\it We can not get the global well-posedness for
the case $0<\alpha\leq 1$, since that in this case one need a
higher order a priori estimate than $L^2$. While in the case
$1<\alpha\leq 3$, we need only the a priori estimate in $L^2$. It
will be a interesting question to find out a suitable a priori
estimate for the dissipative vision of the mKdV and then extend
the local well-posedness to the global well-posedness when the
dissipative term is very weak.}


This paper is organized as follows. In Section 2 we prove linear
estimates in the function space $X^{s,\,1/2}_\alpha$ and in Section
3 we introduce Tao's $[k;\,Z]-$multiplier norm estimate and derive
some trilinear estimate for the nonlinear term $\partial_x(u^3)$
from Tao's estimate. For Section 4, we consider the local
well-posedness while the global well-posedness will be in Section 5.

\section{Linear Estimates}
 In this section we study the linear operator
 $\psi(\cdot)W(\cdot)$ as well as the linear operator $L$ defined
 by
 \begin{eqnarray}
 L:\,f\mapsto\chi_{{\Bbb R}_+}(t)\psi(t)\int_0^tW(t-t')f(t')dt'.
 \end{eqnarray}
 \subsection{Linear estimate for the free term}

\newtheorem{lemma}{Lemma}
\begin{lemma}  Let $s\in {\Bbb R}$. There exists $C>0$ such that
\begin{eqnarray}
\|\psi(t)W(t)\varphi\|_{X^{s,\,1/2}_\alpha}\leq
C\|\varphi\|_{H^s},\,\forall\varphi\in H^s({\Bbb R}).
\end{eqnarray}
\end{lemma}
{\it Proof}.\, By definition of $\|\cdot\|_{X^{s,\,1/2}_\alpha}$,
\begin{eqnarray}
&&\|\psi(t)W(t)\varphi\|_{X^{s,\,1/2}_\alpha}\nonumber\\
&&=\Big\|<\xi>^s<i(\tau-\xi^3)+|\xi|^\alpha>^{1/2}{\cal
F}_t\Big(\psi(t)e^{-|t||\xi|^\alpha}e^{it\xi^3}\hat{\varphi}
(\xi)\Big)(\tau)\Big\|_{L^2_{\xi,\,\tau}}\nonumber\\
&&=\left\|<\xi>^s\hat{\varphi}(\xi)\Big\|<i\tau+|\xi|^\alpha>^{1/2}
{\cal
F}_t\Big(\psi(t)e^{-|t||\xi|^\alpha}\Big)(\tau)\Big\|_{L^2_\tau}\right\|_{L^2_\xi}\nonumber\\
&&\leq C\left\|<\xi>^s\hat{\varphi}(\xi)\Big\|<\tau>^{1/2} {\cal
F}_t\Big(\psi(t)e^{-|t||\xi|^\alpha}\Big)(\tau)\Big\|_{L^2_\tau}\right\|_{L^2_\xi}\nonumber\\
&&\quad + C\left\|<\xi>^{s+\alpha/2}\hat{\varphi}(\xi)\Big\|{\cal
F}_t\Big(\psi(t)e^{-|t||\xi|^\alpha}\Big)(\tau)\Big\|_{L^2_\tau}\right\|_{L^2_\xi}\nonumber\\
&&\leq C\|\varphi\|_{H^s},
\end{eqnarray}
where we used the fact
$$\|\psi(t)e^{-|t||\xi|^\alpha}\|_{H^b_t}\leq C<\xi>^{\frac\alpha
2(2b-1)},\,\forall 0\leq b\leq 1,$$ which can be obtained from
(19) in \cite{MR2}.

\subsection{Linear estimates for the forcing term}
\begin{lemma} For $w\in {\cal S}({\Bbb R}^2)$ consider $k_\xi$
defined on ${\Bbb R}$ by
\begin{eqnarray}
k_\xi(t)=\psi(t)\int_{{\Bbb
R}}\frac{e^{it\tau}-e^{-|t||\xi|^\alpha}}{i\tau+|\xi|^\alpha}\hat{w}(\xi,\,\tau)d\tau.
\end{eqnarray}
Then, it holds for all $\xi\in{\Bbb R}$
\begin{eqnarray}
&&\Big\|<i\tau+|\xi|^\alpha>^{1/2}{\cal
F}_t(k_\xi)\Big\|^2_{L^2({\Bbb R})}\nonumber\\
&&\leq C\Big [\Big(\int_{{\Bbb
R}}\frac{|\hat{w}(\xi,\,\tau)|}{<i\tau+|\xi|^\alpha>}d\tau\Big)^2+\Big(\int_{{\Bbb
R}}\frac{|\hat{w}(\xi,\,\tau)|^2}{<i\tau+|\xi|^\alpha>}d\tau\Big)\Big]
\end{eqnarray}
\end{lemma}
By a little modification of Proposition 2 in \cite{MR2}, we can
obtain the proof of Lemma 2.
\begin{lemma}
Let $s\in{\Bbb R}$.

a)\, There exists $C>0$ such that for all $\upsilon\in {\cal
S}({\Bbb R}^2)$,
\begin{eqnarray}
&&\Big\|\chi_{{\Bbb
R}_+}(t)\psi(t)\int_0^tW(t-t')\upsilon(t')dt'\Big\|_{X^{s,\,1/2}_\alpha}\nonumber\\
&&\quad\quad\leq
C\Big[\|\upsilon\|_{X^{s,\,-1/2}_\alpha}+\Big(\int<\xi>^{2s}(\int\frac{|\hat{\upsilon}
(\xi,\,\tau+\xi^3)|}{<i\tau+|\xi|^\alpha>}d\tau)^2d\xi\Big)^{1/2}\Big]
\end{eqnarray}

b) For any $0<\delta<1/2$ there exists $C_\delta>0$ such that for
all $\upsilon\in X^{s,\,-1/2+\delta}_\alpha$
\begin{eqnarray}
\Big\|\chi_{{\Bbb
R}_+}(t)\psi(t)\int_0^tW(t-t')\upsilon(t')dt'\Big\|_{X^{s,\,1/2}_\alpha}
\leq C_\delta \|\upsilon\|_{X^{s,\,-1/2+\delta}_\alpha}.
\end{eqnarray}
\end{lemma}
{\it Proof}.\, Assume that $\upsilon\in{\cal S}({\Bbb R}^2)$.
Taking the $x-$Fourier transform we get
\begin{eqnarray*}
&&\chi_{{\Bbb
R}_+}(t)\psi(t)\int_0^tW(t-t')\upsilon(t')dt'\\
&&\quad\quad =U(t)\Big[\chi_{{\Bbb R}_+}(t)\psi(t)\int_{{\Bbb
R}}e^{ix\xi}\int^t_0e^{-|t-t'||\xi|^\alpha}{\cal
F}_x(U(-t')\upsilon(t'))(\xi)dt'd\xi\Big].
\end{eqnarray*}
Set $\omega(t')=U(-t')\upsilon(t')$. Writing ${\cal
F}_x(\omega)(\xi,\,t)$ with the help of its time Fourier transform
and using Fubini's theorem one infers that
\begin{eqnarray}
&&\chi_{{\Bbb
R}_+}(t)\psi(t)\int_0^tW(t-t')\upsilon(t')dt'\nonumber\\
&&\quad\quad =U(t)\Big[\chi_{{\Bbb R}_+}(t)\psi(t)\int_{{\Bbb
R}^2}e^{ix\xi}e^{-t|\xi|^\alpha}\hat{\omega}(\xi,\,\tau)\int^t_0e^{it'\tau}
e^{t'|\xi|^\alpha}dt'd\xi d\tau\Big]\nonumber\\
&&\quad\quad =U(t)\Big[\chi_{{\Bbb R}_+}(t)\psi(t)\int_{{\Bbb
R}^2}e^{ix\xi}\frac{e^{it\tau}-e^{-t|\xi|^\alpha}}{i\tau+|\xi|^\alpha}
\hat{\omega}(\xi,\,\tau)d\xi d\tau\Big].
\end{eqnarray}
We set $k_\xi(t)=\psi(t)\int_{\Bbb
R}\frac{e^{it\tau}-e^{-|t||\xi|^\alpha}}{i\tau+|\xi|^\alpha}\hat{\omega}
(\xi,\,\tau)d\tau$. Since $\omega(t)=U(-t)\upsilon(t)\in{\cal
S}({\Bbb R}^2)$, it is clear that for any fixed $\xi\in{\Bbb
R},\,k_\xi$ is continuous on ${\Bbb R}$ and $k_\xi(0)=0$. Then it
is not too hard to derive that $\|\chi_{{\Bbb
R}_+}(t)k_\xi(t)\|_{H^b_t}\leq \|k_\xi(t)\|_{H^b_t},\,0\leq b\leq
1$, and $b\neq 1/2$. The case $b=1/2$ follows by Lebesgue
dominated convergence theorem. Thus in view of (20)
\begin{eqnarray*}
&&\Big\|\chi_{{\Bbb
R}_+}(t)\psi(t)\int_0^tW(t-t')\upsilon(t')dt'\Big\|_{X^{s,\,1/2}_\alpha}\\
&&=\Big\|<i\tau+|\xi|^\alpha>^{1/2}<\xi>^s{\cal F}_{x,\,t}
\Big(U(-t)\chi_{{\Bbb
R}_+}(t)\psi(t)\int_0^tW(t-t')\upsilon(t')dt'\Big)\Big\|_{L^2_{\xi,\,\tau}}\\
&&\leq \left\|<\xi>^s\Big\|{\cal F}_{x} \Big(U(-t)\chi_{{\Bbb
R}_+}(t)\psi(t)\int_0^tW(t-t')\upsilon(t')dt'\Big)\Big\|_{H_\tau^{1/2}}\right\|_{L^2_{\xi}}\\
&&\quad +\left\|<\xi>^{s+\alpha/2}\Big\|{\cal F}_{x}
\Big(U(-t)\chi_{{\Bbb
R}_+}(t)\psi(t)\int_0^tW(t-t')\upsilon(t')dt'\Big)\Big\|_{L_\tau^{2}}\right\|_{L^2_{\xi}}\\
&&\leq \Big\|<\xi>^s\|\chi_{{\Bbb
R}_+}(t)k_\xi(t)\|_{H^{1/2}_t}\Big\|_{L^2_\xi}+\Big\|<\xi>^{s+\alpha/2}\|\chi_{{\Bbb
R}_+}(t)k_\xi(t)\|_{L^{2}_t}\Big\|_{L^2_\xi}\\
&&\leq
\Big\|<\xi>^s\|k_\xi(t)\|_{H^{1/2}_t}\Big\|_{L^2_\xi}+\Big\|<\xi>^{s+\alpha/2}\|
k_\xi(t)\|_{L^{2}_t}\Big\|_{L^2_\xi}\\
&&\leq C\left\|<\xi>^s\Big\|<i\tau+|\xi|^\alpha>^{1/2}{\cal
F}_t(k_\xi(t)\Big\|_{L^2_\tau}\right\|_{L^2_\xi}.
\end{eqnarray*}
Then (18) follows directly from Lemma 2 together with the last
estimate.\\
To prove (19) we first assume that $\upsilon\in{\cal S}({\Bbb
R}^2)$. Applying Cauchy-Schwartz inequality in $\tau$ on the
second term of the right-hand side of (18), one obtains (19) for
$\upsilon\in{\cal S}({\Bbb R}^2)$. The result for $\upsilon\in
X^{s,\,-1/2+\delta}_\alpha$ follows by density.

\section{Tao's $[k;\,Z]-$multiplier norm estimate and its application}

In this section we introduce Tao's $[k;\,Z]-$multiplier norm
estimate and derive the trilinear estimate needed to obtain the
local existence result from Tao's multiplier norm estimate for KdV
equation.

Let $Z$ be any abelian additive group with an invariant measure
$d\xi$. For any integer $k\geq 2$, we let $\Gamma_k(Z)$ denote the
hyperplane
$$\Gamma_k(Z):=\{(\xi_1,\cdots,\,\xi_k)\in
Z^k:\,\xi_1+\cdots+\xi_k=0\}$$ which is endowed with the measure
$$\int_{\Gamma_k(Z)}f:=\int_{Z^{k-1}}f(\xi_1,\cdots,\,\xi_{k-1},\,-\xi_1
-\cdots-\xi_{k-1})d\xi_1\cdots d\xi_{k-1}.$$

A $[k;\,Z]-$multiplier is defined to be any function
$m:\,\Gamma_k(Z)\rightarrow {\Bbb C}$ which was introduced by Tao
in \cite{Tao}. And the multiplier norm $\|m\|_{[k;\,Z]}$ is
defined to be the best constant such that the inequality
\begin{eqnarray}
\Big |\int_{\Gamma_k(Z)}m(\xi)\prod_{j=1}^kf_j(\xi_j)\Big |\leq
\|m\|_{[k;\,Z]}\prod_{j=1}^k\|f_j\|_{L^2(Z)},
\end{eqnarray}
holds for all test functions $f_j$ on $Z$. Tao systematically
studied this kind of weighted convolution estimates on $L^2$ in
\cite{Tao}. To state Tao's results, we use some notation he used
in his paper.

We use $A\lesssim B$ to denote the statement that $A\leq CB$ for
some large constant $C$ which may vary from line to line and
depend on various parameters, and similarly use $A\ll B$ to denote
the statement $A\leq C^{-1}B$. We use $A\sim B$ to denote the
statement that $A\lesssim B\lesssim A$.

Any summations over capitalized variables such as $N_j,\,L_j,\,H$
are presumed to be dyadic, i.e., these variables range over numbers
of the form $2^k$ for $k\in {\Bbb Z}$. In this paper, we will only
consider the $[3;\,Z]$-multiplier. Let $N_1,\,N_2,\,N_3>0$. It will
be convenient to define the quantities $N_{max}\geq N_{med}\geq
N_{min}$ to be the maximum, median, and minimum of $N_1,\,N_2,\,N_3$
respectively. Similarly define $L_{max}\geq L_{med}\geq L_{min}$
whenever $L_1,\,L_2,\,L_3>0$. And we also adopt the following
summation conventions. Any summation of the form $L_{max}\sim\cdots$
is a sum over the three dyadic variables $L_1,\,L_2,\,L_3\gtrsim 1$,
thus for instance
$$\sum_{L_{max}\sim H}:=\sum_{L_1,\,L_2,\,L_3\gtrsim 1:\,L_{max}\sim
H}.$$ Similarly, any summation of the form $N_{max}\sim\cdots$ sum
over the three dyadic variables $N_1,\,N_2,\,N_3>0$, thus for
instance
$$\sum_{N_{max}\sim N_{med}\sim
N}:=\sum_{N_1,\,N_2,\,N_3>0:\,N_{max}\sim N_{med}\sim N}.$$ If
$\tau,\,\xi$ and $h(\cdot)$ are given, we also adopt the convention
that $\lambda$ is short-hand for $$\lambda:=\tau-h(\xi).$$ Similarly
we have \begin{eqnarray} \lambda_j:=\tau_j-h_j(\xi_j).
\end{eqnarray}

In this paper, we do not go further on the general framework of
Tao's weighted convolution estimates. We focus our attention to the
$[k;\,Z]-$multiplier norm estimate for KdV equation. During the
estimate we need the resonance function
\begin{eqnarray}
h(\xi)=\xi_1^3+\xi_2^3+\xi_3^3=-\lambda_1-\lambda_2-\lambda_3,
\end{eqnarray}
which measures to what extent the spatial frequencies
$\xi_1,\,\xi_2,\,\xi_3$ can resonate with each other.

By dyadic decomposition of the variables $\xi_j,\,\lambda_j$, as
well as the function $h(\xi)$, one is led to consider
\begin{eqnarray}
\|X_{N_1,\,N_2,\,N_3;\,H;\,L_1,\,L_2,\,L_3}\|_{[3,\,{\Bbb
R}\times{\Bbb R}]},
\end{eqnarray}
where $X_{N_1,\,N_2,\,N_3;\,H;\,L_1,\,L_2,\,L_3}$ is the
multiplier
\begin{eqnarray}
X_{N_1,\,N_2,\,N_3;\,H;\,L_1,\,L_2,\,L_3}(\xi,\,\tau):=\chi_{|h(\xi)|\sim
H}\prod_{j=1}^3\chi_{|\xi_j|\sim N_j}\chi_{|\lambda_j|\sim L_j}.
\end{eqnarray}

From the identities
$$\xi_1+\xi_2+\xi_3=0$$
and
$$\lambda_1+\lambda_2+\lambda_3+h(\xi)=0$$
on the support of the multiplier, we see that
$X_{N_1,\,N_2,\,N_3;\,H;\,L_1,\,L_2,\,L_3}$ vanishes unless
\begin{eqnarray}
N_{max}\sim N_{med},
\end{eqnarray}
and
\begin{eqnarray}
L_{max}\sim \max(H,\,L_{med}).
\end{eqnarray}
From the resonance identity
\begin{eqnarray}
h(\xi)=\xi_1^3+\xi_2^3+\xi_3^3=3\xi_1\xi_2\xi_3
\end{eqnarray}
we see that we may assume that
\begin{eqnarray}
H\sim N_1N_2N_3,
\end{eqnarray}
since the multiplier in (24) vanishes otherwise.

Now we are in the position to state Tao's $[k;\,Z]-$multiplier
norm estimate for KdV equation in non-periodic case.
\begin{lemma}(see Proposition 6.1 in \cite{Tao}).\, Let
$H,\,N_1,\,N_2,\,N_3,\,L_1,\,L_2,\,L_3>0$ obey (26), (27), (29).

$\bullet$((++)Coherence) If $N_{max}\sim N_{min}$ and $L_{max}\sim
H$, then we have
\begin{eqnarray}
(24)\lesssim L_{min}^{1/2}N_{max}^{-1/4}L_{med}^{1/4}.
\end{eqnarray}

$\bullet$((+-)Coherence) If $N_2\sim N_3\gg N_1$ and $H\sim
L_1\gtrsim L_2,\,L_3$, then
\begin{eqnarray}
(24)\lesssim
L_{min}^{1/2}N_{max}^{-1}\min(H,\,\frac{N_{max}}{N_{min}}L_{med})^{1/2}.
\end{eqnarray}
Similarly for permutations.

$\bullet$ In all other cases, we have
\begin{eqnarray}
(24)\lesssim L_{min}^{1/2}N_{max}^{-1}\min(H,\,L_{med})^{1/2}.
\end{eqnarray}
\end{lemma}

With Lemma 4 one can now derive a trilinear estimate involving the
new Bourgain spaces $X^{s,\,b}_\alpha$.

\begin{lemma} \, Let $s>1/4-\alpha/4$. For all $u_1,\,u_2,\,u_3$ on ${\Bbb
R}\times {\Bbb R}$ and $0<\epsilon\ll 1$, we have
\begin{eqnarray}
\|\partial_x(u_1u_2u_3)\|_{X^{s,\,-1/2+\epsilon}_\alpha}\lesssim
\|u_1\|_{X^{s,\,1/2}_\alpha}\|u_2\|_{X^{s,\,1/2}_\alpha}\|u_3\|_{X^{s,\,1/2}_\alpha}.
\end{eqnarray}
\end{lemma}

As seen in \cite{MR2}, Theorem 1 can be reduced to the trilinear
estimate Lemma 5 and the linear estimates we obtained in section
2. What we need to do in the following is to prove the trilinear
estimate.

{\it Proof of Lemma 5.}\, By duality and Plancherel it suffices to
show that
\begin{eqnarray*}
\left\|\frac{(\xi_1+\xi_2+\xi_3)<\xi_4>^s}{<\tau_4-\xi_4^3+i|\xi_4|^\alpha>^{1/2-\epsilon}
\prod_{j=1}^3<\xi_j>^s<\tau_j-\xi_j^3+i|\xi_j|^\alpha>^{1/2}}\right\|_{[4,\,{\Bbb
R}\times{\Bbb R}]}\lesssim 1.
\end{eqnarray*}
We estimate $|\xi_1+\xi_2+\xi_3|$ by $<\xi_4>$. We then apply the
inequality fractional Leibnitz rule
$$<\xi_4>^{s+1}\lesssim <\xi_4>^{1/2}\sum_{j=1}^3<\xi_j>^{s+1/2}$$
where we assume $-1/2<s<1/2$, and symmetry to reduce to
\begin{eqnarray*}
\left\|\frac{<\xi_1>^{-s}<\xi_3>^{-s}<\xi_2>^{1/2}<\xi_4>^{1/2}}{<\tau_4-\xi_4^3+i|\xi_4|^\alpha>^{1/2-\epsilon}
\prod_{j=1}^3<\tau_j-\xi_j^3+i|\xi_j|^\alpha>^{1/2}}\right\|_{[4,\,{\Bbb
R}\times{\Bbb R}]}\lesssim 1.
\end{eqnarray*}
We may replace $<\tau_2-\xi_2^3+i|\xi_2|^\alpha>^{1/2}$ by
$<\tau_2-\xi_2^3+i|\xi_2|^\alpha>^{1/2-\epsilon}$. By the $TT^*$
identity (see Lemma 3.7, p847 in\cite{Tao}), the estimate is
reduced to the following bilinear estimate.

\begin{lemma}(Bilinear estimate).\, Let $s>1/4-\alpha/4$. For all $u,\,v$ on ${\Bbb
R}\times {\Bbb R}$ and $0<\epsilon\ll 1$, we have
\begin{eqnarray}
\|uv\|_{L^2({\Bbb R}\times{\Bbb R})}\lesssim
\|u\|_{X^{-1/2,\,1/2-\epsilon}_\alpha({\Bbb R}\times{\Bbb R})
}\|v\|_{X^{s,\,1/2}_\alpha({\Bbb R}\times{\Bbb R})}.
\end{eqnarray}
\end{lemma}

{\it Proof of the bilinear estimate}.\, By Plancherel it suffices
to show that
\begin{eqnarray}
\left\|\frac{<\xi_1>^{-s}<\xi_2>^{1/2}}{<\tau_1-\xi_1^3+i|\xi_1|^\alpha>^{1/2}
<\tau_2-\xi_2^3+i|\xi_2|^\alpha>^{1/2-\epsilon}}\right\|_{[3,\,{\Bbb
R}\times{\Bbb R}]}\lesssim 1.
\end{eqnarray}

There are two approaches to computing these quantities and their
generalization. One approach proceeds using the Cauchy-Schwarz
inequality, this reducing matters to integrating certain weights
on intersections of hypersurfaces $\tau=h(\xi)$; the other
utilizes dyadic decomposition and orthogonality before resorting
to Cauchy-Schwarz. The advantages of dyadic decomposition are that
one can re-use the estimates on dyadic blocks to prove other
estimates, and the nature of interactions between different scales
of frequency is more apparent.

By dyadic decomposition of the variables
$\xi_j,\,\lambda_j,\,h(\xi)$, we may assume that $|\xi_j|\sim
N_j,\,|\lambda_j|\sim L_j,\,|h(\xi)|\sim H$. By the translation
invariance of the $[k;\,Z]$-multiplier norm, we can always restrict
our estimate on $\lambda_j\gtrsim 1$ and $\max(N_1,N_2,N_3)\gtrsim
1$. The comparison principle and orthogonality (see Schur's test in
\cite{Tao} p851) reduce our estimate to show that
\begin{eqnarray}
&&\sum_{N_{max}\sim N_{med}\sim N}\sum_{L_1,\,L_2,\,L_3\gtrsim
1}\frac{<N_1>^{-s}<N_2>^{1/2}}{\max(L_1,\,<N_1>^\alpha)^{1/2}
\max(L_2,\,<N_2>^\alpha)^{1/2-\epsilon}}\nonumber\\
&&\quad\qquad\quad\quad\quad\qquad\left\|X_{N_1,\,N_2,\,N_3;\,L_{max};\,L_1,\,L_2,\,L_3}\right\|_{[3;\,{\Bbb
R}\times{\Bbb R}]}\lesssim 1
\end{eqnarray}
and
\begin{eqnarray}
&&\sum_{N_{max}\sim N_{med}\sim N}\sum_{L_{max}\sim
L_{med}}\sum_{H\ll L_{max}}
\frac{<N_1>^{-s}<N_2>^{1/2}}{\max(L_1,\,<N_1>^\alpha)^{1/2}
\max(L_2,\,<N_2>^\alpha)^{1/2-\epsilon}}\nonumber\\
&&\quad\qquad\quad\quad\quad\qquad\left\|X_{N_1,\,N_2,\,N_3;\,H;\,L_1,\,L_2,\,L_3}\right\|_{[3;\,{\Bbb
R}\times{\Bbb R}]}\lesssim 1
\end{eqnarray}
for all $N\gtrsim 1$. This can be accomplished by Tao's estimate
on dyadic blocks for KdV equation, i.e., Lemma 4 and some concrete
summation.

Fix $N\gtrsim 1$. We first prove (37). We may assume (29). By (32)
we reduce to
\begin{eqnarray*}
&&\sum_{N_{max}\sim N_{med}\sim N}\sum_{L_{max}\sim L_{med}\gtrsim
N_1N_2N_3}\\
&&\quad\quad\frac{<N_1>^{-s}<N_2>^{1/2}}{\max(L_1,\,<N_1>^\alpha)^{1/2}
\max(L_2,\,<N_2>^\alpha)^{1/2-\epsilon}}
L_{min}^{1/2}N_{min}^{1/2}\lesssim 1.
\end{eqnarray*}
This high modulation case is easier to handle. We do not need to
consider the effect of the dissipative term, i.e., the part with
symbol $|\xi|^\alpha$. Just crudely estimating
$$
<N_1>^{-s}<N_2>^{1/2}\lesssim N^{1/2+\max(0,\,-s)}$$ and
$$
\max(L_1,\,<N_1>^\alpha)^{1/2}\max(L_2,\,<N_2>^\alpha)^{1/2-\epsilon}
\gtrsim L_{min}^{1/2}L_{med}^{1/2-\epsilon},
$$
then performing the $L$ summations, we reduce to
\begin{eqnarray*}
\sum_{N_{max}\sim N_{med}\sim
N}\frac{N_{min}^\epsilon}{N^{1/2-\max(0,\,-s)-2\epsilon}}\lesssim
1,
\end{eqnarray*}
which is true if $1/2-\max(0,\,-s)>0$. So, (37) is true if
$s>-1/2$.

Now we show the low modulation case (36). We may assume
$L_{max}\sim N_1N_2N_3$. We first deal with the contribution where
(30) holds. In this case we have $N_1,\,N_2,\,N_3\sim N\gtrsim 1$,
so we reduce to
\begin{eqnarray}
\sum_{L_{max}\sim
N^3}\frac{N^{-s}N^{1/2}}{\max(L_1,\,N^\alpha)^{1/2}\max(L_2,\,N^\alpha)^{1/2-\epsilon}}
L_{min}^{1/2}N^{-1/4}L_{med}^{1/4}\lesssim 1.
\end{eqnarray}
Note that
$$\max(L_1,\,N^\alpha)^{1/2}\max(L_2,\,N^\alpha)^{1/2-\epsilon}\gtrsim
\max(L_{min},\,N^\alpha)^{1/2}\max(L_{med},\,N^\alpha)^{1/2-\epsilon}.$$
We may assume that $L_1\leq L_2\leq L_3$. We now need consider
three subcases: $L_2\leq N^\alpha,\,L_1\leq N^\alpha\leq
L_2,\,N^\alpha\leq L_1$.

If $L_2\leq N^\alpha$, we reduce to
\begin{eqnarray*}
\sum_{L_1\leq L_2\leq N^\alpha\lesssim L_3\sim N^3
}\frac{N^{-s}N^{1/4}L_1^{1/2}L_2^{1/4}}{N^{\alpha/2}N^{\alpha(1/2-\epsilon)}}
\lesssim\frac{N^{1/4-s}N^{3\alpha/4}}{N^{\alpha-\alpha\epsilon}}\lesssim
1,
\end{eqnarray*}
which is true if $1/4-s-\alpha/4<0$, i.e., $s>1/4-\alpha/4$.

If $L_1\leq N^\alpha\leq L_2$, we reduce to
\begin{eqnarray*}
\sum_{L_1\leq N^\alpha\leq L_2\leq L_3\sim N^3
}\frac{N^{-s}N^{1/4}L_1^{1/2}L_2^{1/4}}{N^{\alpha/2}L_2^{1/2-\epsilon}}
\lesssim\frac{N^{1/4-s}N^{\alpha/2}}{N^{\alpha/2}N^{(1/4-\epsilon)\alpha}}\lesssim
1,
\end{eqnarray*}
which holds if $s>1/4-\alpha/4$.

If $N^\alpha\leq L_1$, we reduce to
\begin{eqnarray*}
\sum_{N^\alpha\leq L_1\leq L_2\leq L_3\sim N^3
}\frac{N^{-s}N^{1/4}L_1^{1/2}L_2^{1/4}}{L_1^{1/2}L_2^{1/2-\epsilon}}
\lesssim\frac{N^{1/4-s}}{N^{(1/4-\epsilon)\alpha}}\lesssim 1,
\end{eqnarray*}
which is true if $s>1/4-\alpha/4$.

Now we deal with the cases where (31) applies. We do not have
perfect symmetry and must consider three cases
$$\begin{array}{ll}
N\sim N_1\sim N_2\gg N_3;\,& H\sim L_3\gtrsim L_1,\,L_2\\
N\sim N_2\sim N_3\gg N_1;\,& H\sim L_1\gtrsim L_2,\,L_3\\
N\sim N_1\sim N_3\gg N_2;\,& H\sim L_2\gtrsim L_1,\,L_3
\end{array}$$
separately.

In the first case we reduce by (31) to
\begin{eqnarray*}
&&\sum_{N_3\ll N}\sum_{1\lesssim L_1,\,L_2\lesssim
N^2N_3}\frac{N^{1/2-s}}{\max(L_1,\,N^\alpha)^{1/2}\max(L_2,\,N^\alpha)^{1/2-\epsilon}}\\
&&\quad\qquad\quad\quad\qquad
L_{min}^{1/2}N^{-1}\min(N^2N_3,\,\frac{N}{N_3}L_{med})^{1/2}\lesssim
1.
\end{eqnarray*}
Performing the $N_3$ summation we  reduce to
\begin{eqnarray*}
\sum_{1\lesssim L_1,\,L_2\lesssim
N^3}\frac{N^{1/2-s}}{\max(L_1,\,N^\alpha)^{1/2}\max(L_2,\,N^\alpha)^{1/2-\epsilon}}
L_{min}^{1/2}N^{-1}N^{3/4}L_{med}^{1/4}\lesssim 1.
\end{eqnarray*}
which is the same as (38) and so it is true if $s>1/4-\alpha/4$.

We can deal with the second and third cases in a unified way. By
asymmetry it suffices to show the worst case. We simplify using
the first half of (31) to
$$\sum_{N_{min}\ll N}\sum_{1\lesssim L_{min},\,L_{med}\ll
N^2N_{min}}\frac{<N_{min}>^{-s}N^{1/2}}{L_{min}^{1/2}L_{max}^{1/2-\epsilon}}
L_{min}^{1/2}N_{min}^{1/2}\lesssim 1.$$ We may assume
$N_{min}\gtrsim N^{-2}$ since the inner sum vanishes otherwise.
Performing the $L$ summation we reduce to
$$\sum_{N^{-2}\lesssim N_{min}\ll
N}\frac{<N_{min}>^{-s}N^{1/2}N_{min}^{1/2}}{(N^2N_{min})^{1/2-\epsilon}}\lesssim
1$$ which holds if $s>-1/2$.

To finish the proof of (36) it remains to deal with the cases
where (32) holds. This reduces to
\begin{eqnarray*}
&&\sum_{N_{max}\sim N_{med}\sim N}\sum_{L_{max}\sim N_1N_2N_3}\\
&&\quad\quad\frac{<N_1>^{-s}<N_2>^{1/2}}{\max(L_1,\,<N_1>^\alpha)^{1/2}
\max(L_2,\,<N_2>^\alpha)^{1/2-\epsilon}}
L_{min}^{1/2}N^{-1}L_{med}^{1/2}\lesssim 1.
\end{eqnarray*}
Performing the $L$ summations, we reduce to
$$\sum_{N_{max}\sim N_{med}\sim
N}<N_1>^{-s}<N_2>^{1/2}(N_1N_2N_3)^\epsilon N^{-1}\lesssim 1$$
which is true if $s>-1/2$.

\section{Local well-posedness}
 Let $\varphi$ in $H^s({\Bbb
R}),\,s>\frac{1}{4}-\frac{\alpha}{4}$. We first prove the existence
of a solution $u$ of the integral formulation (10) of the equation
(2) on some interval $[0,\,T]$ for $T\leq 1$ small enough. Clearly,
if $u$ is a solution of the integral equation $u=F(u)$ with
\begin{eqnarray}
F(u)=\psi(t)\Big[W(t)\varphi-\chi_{{\Bbb R}_+}(t)
\int_0^tW(t-t')\partial_x(\psi_T^3(t')u^3(t'))dt'\Big],
\end{eqnarray}
then $u$ is a solution of (10) on $[0,\,T]$. We are going to solve
(39) in the space
\begin{eqnarray}
Z=\{u\in
X^{s,\,1/2},\,\|u\|_Z=\|u\|_{X^{s_c^+,\,1/2}}+\nu\|u\|_{X^{s,\,1/2}}<+\infty\},
\end{eqnarray}
where $s_c^+\in(\frac{1}{4}-\frac{\alpha}{4},\,s)$ is fixed, and
where the constant $\nu$ is defined for all nontrivial $\varphi$ by
\begin{eqnarray}
\nu=\frac{\|\varphi\|_{H^{s_c^+}}}{\|\varphi\|_{H^s}}.
\end{eqnarray}
To get an contractive factor in our argument, we need the following
modified trilinear estimates.
\begin{lemma}\label{mainlema} \, Given $s>\frac{1}{4}-\frac{\alpha}{4}$, there exist
$C,\,\mu,\,\delta>0$ such that for any triple $(u,\,v,\,w)\in
X^{s,\,1/2}$ with compact support in $[-T,\,T]$
\begin{eqnarray}\label{triest}
\|\partial_x(uvw)\|_{X^{s,\,-1/2+\delta}}\leq CT^\mu
\|u\|_{X^{s,\,1/2}}\|v\|_{X^{s,\,1/2}}\|w\|_{X^{s,\,1/2}}.
\end{eqnarray}
\end{lemma}
The following lemma is a direct consequence of Lemma \ref{mainlema}
together with the triangle inequality
\begin{eqnarray}
\forall s\geq s_c^+,\,\,<\xi>^s&\leq &
<\xi>^{s_c^+}<\xi_1>^{s-s_c^+}+<\xi>^{s_c^+}<\xi_2>^{s-s_c^+}\nonumber\\&&\quad+
<\xi>^{s_c^+}<\xi-\xi_1-\xi_2>^{s-s_c^+}.
\end{eqnarray}

\begin{lemma} \label{uselema}\, Given $s_c^+>\frac{1}{4}-\frac{\alpha}{4}$,
there exist $C,\,\mu,\,\delta>0$ such that for any $s\geq s_c^+$ and
any triple $(u,\,v,\,w)\in X^{s,\,1/2}$ with compact support in
$[-T,\,T]$
\begin{eqnarray}\label{mtrist}
&&\|\partial_x(uvw)\|_{X^{s,\,-1/2+\delta}}\leq CT^\mu
\Big(\|u\|_{X^{s,\,1/2}}\|v\|_{X^{s_c^+,\,1/2}}\|w\|_{X^{s_c^+,\,1/2}}+\|u\|_{X^{s_c^+,\,1/2}}
\nonumber\\
&&\quad\|v\|_{X^{s,\,1/2}}\|w\|_{X^{s_c^+,\,1/2}}
+\|u\|_{X^{s_c^+,\,1/2}}\|v\|_{X^{s_c^+,\,1/2}}\|w\|_{X^{s,\,1/2}}\Big).
\end{eqnarray}
\end{lemma}

The proof of Lemma \ref{mainlema} is exactly the same with that of
Lemma 5 except changing $f_4(\xi_4,\,\tau_4)$ into
$\frac{f_4(\xi_4,\,\tau_4)}{<|\tau_4-\xi_4^3|+|\xi_4|^{\alpha}>^{\delta}}$
and using the following lemma due to Molinet and Ribaud \cite{MR2}.
 \begin{lemma}
 \label{Lemma 4} Let $f$ with support in $[-T,\,T]$ in time. For
 any $\delta>0$, there exists $\mu=\mu(\delta)>0$ such that
 $$\Big\|\frac{\hat{f}(\xi,\,\tau)}{<\tau-\xi^3>^\delta}\Big\|_{L^2_{\xi,\,\tau}}\leq
 C T^\mu\|f\|_{L^2_{t,\,x}}.$$
 \end{lemma}

By Lemmas 7 and 8, there exist $\delta,\,\mu>0$ only depending on
$s_c^+$ such that
\begin{eqnarray}
&&\|F(u)\|_{X^{s_c^+,\,1/2}}\leq
C\|\varphi\|_{H^{s_c^+}}+CT^\mu\|u\|_{X^{s_c^+,\,1/2}}^3,\\
&&\|F(u)\|_{X^{s,\,1/2}}\leq
C\|\varphi\|_{H^{s}}+CT^\mu\|u\|_{X^{s_c^+,\,1/2}}^2\|u\|_{X^{s,\,1/2}}.
\end{eqnarray}
Gathering the above two estimates, one infers that
\begin{eqnarray}
\|F(u)\|_Z\leq
C(\|\varphi\|_{H^{s_c^+}}+\nu\|\varphi\|_{H^s})+CT^\mu\|u\|_Z^3.
\end{eqnarray}
Next, since
$\partial_x(u^3)-\partial_x(v^3)=\partial_x[(u-v)(u^2+uv+v^2)]$, we
get in the same way that
\begin{eqnarray}
\|F(u)-F(v)\|_{X^{s_c^+,\,1/2}}&\leq &
CT^\mu\|u-v\|_{X^{s_c^+,\,1/2}}
\Big(\|u\|_{X^{s_c^+,\,1/2}}^2\nonumber\\&&\quad
+\|u\|_{X^{s_c^+,\,1/2}}
\|v\|_{X^{s_c^+,\,1/2}}+\|v\|_{X^{s_c^+,\,1/2}}^2\Big),
\end{eqnarray}
and \begin{eqnarray}
 \|F(u)-F(v)\|_{X^{s,\,1/2}}&\leq &
CT^\mu\|u-v\|_{X^{s_c^+,\,1/2}}\Big(\|u\|_{X^{s_c^+,\,1/2}}
\|u\|_{X^{s,\,1/2}}\nonumber\\&&+\|u\|_{X^{s_c^+,\,1/2}}
\|v\|_{X^{s,\,1/2}}+\|u\|_{X^{s,\,1/2}}\|v\|_{X^{s_c^+,\,1/2}}\nonumber\\
&&+
\|v\|_{X^{s_c^+,\,1/2}}\|v\|_{X^{s,\,1/2}}\Big)+CT^\mu\|u-v\|_{X^{s,\,1/2}}
\Big(\|u\|_{X^{s_c^+,\,1/2}}^2\nonumber\\&& +\|u\|_{X^{s_c^+,\,1/2}}
\|v\|_{X^{s^+,\,1/2}}+ \|v\|_{X^{s_c^+,\,1/2}}^2\Big)
\end{eqnarray}
Combining (48) and (49), we deduce that
\begin{eqnarray}
\|F(u)-F(v)\|_Z\leq CT^\mu\|u-v\|_Z\Big(\|u\|_Z+\|v\|_Z\Big)^2.
\end{eqnarray}
Setting $T=\Big(128C^3\|\varphi\|_{H^{s_c^+}}^2\Big)^{-1/\mu}$, we
deduce that from (47) and (50) that $F$ is strictly contractive on
the ball of radius $4C\|\varphi\|_{H^{s_c^+}}$ in $Z$. This proves
the existence of a solution $u\in X^{s,\,1/2}$ to the equation (2)
on the time interval $[0,\,T]$ with
$T=T(\|\varphi\|_{H^{s_c^+}})>0$.

The uniqueness of the solution is exactly the same as in \cite{MR2}.
For brevity, we omit the details.

\vspace{0.5cm}

\section{Global well-posedness for the strong dissipative term.}

We now turn to the global well-posedness for the case $1<\alpha\leq
3$. In this case, we fix
$s_c^+\in(\frac{1}{4}-\frac{\alpha}{4},\,0).$ By an similar argument
as the Proposition 2.4 in \cite{MR2}, we can get the following
lemma.
\begin{lemma}\label{forglobal}Let $s\in\Bbb R$ and $\delta>0$. For
all $f\in X^{s,-1/2+\delta}$,
\begin{equation}
t\longmapsto\int_0^tW(t-t')f(t')dt'\in C(\Bbb R_+,
H^{s+\alpha\delta}).
\end{equation}
Moreover, if $(f_n)$ is a sequence with $(f_n)$ is a sequence with
$f_n\xrightarrow[n\rightarrow\infty]{} 0$ in $X^{s,-1/2+\delta}$,
then
\begin{equation}
\Big\|\int_0^tW(t-t')f_n(t')dt'\Big\|_{L^\infty(\Bbb
R_+,H^{s+\alpha\delta})}\xrightarrow[n\rightarrow\infty]{}0.
\end{equation}
\end{lemma}

It is easily to check that $W(\cdot)\varphi\in
C([0,+\infty);H^s)\cap C((0,+\infty);H^\infty).$ Then it follows
from Lemma 5 ,Lemma 10 and the local existence of the solution that
$$u\in C([0,T]; H^s)\cap C((0;T];H^{s+\alpha\delta}), \text{ for some } T=T(\|\varphi\|_{H^{s^+_c}}).$$
By induction we have $u\in C((0,T];H^\infty)$. Since that the exist
time  $T$ of the solution depends only on the norm of date
$\|\varphi\|_{H^{s^+_c}},$ to extent the solution into the global
sense we need a control on $\|u(t)\|_{H^{s^+_c}}.$ Taking the
$L^2$-scalar product of (2) with $u$, we can get that
$\|u(t)\|_{L^2}$ is nonincreasing on $(0,T]$. This implies that the
solution is global in time.

\hspace{5mm}
\begin{center}
\end{center}


\begin{thebibliography}{99}
 \bibitem{B} D.\,Bekiranov, \it The initial-value problem for the generalized
 Burgers' equation, \rm Differential Integral Equations \bf 9 \rm (1996), 1253-1265.
\bibitem{Bourgain} J.\,Bourgain, \it Fourier transform restriction phenomena
for certain lattice subsets and applications to nonlinear
evolution equations. II. The KdV equation, \rm Geom. Funct. Anal.
 \bf 3 \rm (1993), 209-262.
 \bibitem{CL} W.\, G.\, Chen and J.\, F.\, Li \it On the low
 regularity of the modified Korteweg de Vries equation with a
 dissipative term, \rm submitted.
\bibitem{Dix2} D.\,B.\,Dix, \it Nonuniqueness and uniqueness in the initial-value
problem for Burgers' equation, \rm SIAM J. Math. Anal. \bf 27 \rm
(1996), 708-724.
 \bibitem{KPV1}C.\,E.\,Kenig, G.\,Ponce, and L.\,Vega, \it A bilinear estimate with applications
 to the KdV equation, \rm J. Amer. Math. Soc. \bf 9 \rm (1996), 573-603.
\bibitem{MR1}L.\,Molinet and F.\,Ribaud, \it The Cauchy problem for dissipative Kortewig
de Vries equations in Sobolev spaces of negative order, \rm
Indiana Univ. Math. J. \bf 50 \rm (2001), 1745-1776.
\bibitem{MR2} L.\,Molinet and F.\,Ribaud, \it On the low regularity of the Kortewig de Vries
-Burgers equation, \rm Inter. Math. Research Notices \bf 37 \rm
(2002), 1979-2005.
\bibitem{OS} E.\,Ott and N.Sudan, \it Damping of solitary waves, \rm Phys. Fluids \bf
13 \rm (1970), 1432-1434.
\bibitem{Tao} T.\,Tao, \it Multilinear weighted convolution of
$L^2$ functions, and applications to nonlinear dispersive
equations, \rm Amer. J. of Math. \bf 123 \rm (2001), 839-908.
\end{thebibliography}
\end{document}